\documentclass[12pt,reqno]{amsart}

\usepackage[text={420pt,660pt},centering]{geometry}

\usepackage{amssymb,amsfonts,amsthm,mathrsfs,mathtools}
\usepackage{marginnote}
\usepackage{comment}
\usepackage{enumitem}
\usepackage{graphicx}
\usepackage{bm}
\usepackage{xfrac}
\usepackage{dsfont}
\usepackage[dvipsnames]{color}

\usepackage[colorlinks=true, pdfstartview=FitV, linkcolor=black, citecolor=blue, urlcolor=blue]{hyperref}

\DeclareMathAlphabet{\mathcal}{OMS}{cmsy}{m}{n}

\theoremstyle{plain}

\newtheorem{theorem}{Theorem}[section]

\newtheorem{lemma}[theorem]{Lemma}

\theoremstyle{definition}
\newtheorem{remark}[theorem]{Remark}

\newcommand{\R}{\mathbb{R}}

\newcommand{\Z}{\mathbb{Z}}
\newcommand{\T}{\mathbb{T}}

\newcommand{\p}{\partial}

\let\div\relax
\DeclareMathOperator{\div}{div}

\renewcommand{\T}{\mathbb{T}}

\renewcommand{\d}{\mathrm{d}}

\numberwithin{equation}{section}
\setlist[enumerate]{leftmargin=*}

\title[]{Sharp bounds on enstrophy growth for viscous scalar conservation laws}
\author[D. Albritton]{Dallas Albritton} 
\address[D. Albritton]{University of Wisconsin-Madison, Department of Mathematics,  480 Lincoln Dr, Madison, WI 53706, USA}
\email[]{dalbritton@wisc.edu}

\author[N. De Nitti]{Nicola De Nitti}
\address[N. De Nitti]{Friedrich-Alexander-Universit\"at Erlangen-N\"urnberg, Department of Mathematics, Chair for Dynamics, Control, Machine Learning and Numerics (Alexander von Humboldt Professorship), Cauerstr. 11, 91058 Erlangen, Germany}
\email[]{nicola.de.nitti@fau.de}


\begin{document}

\begin{abstract}
We prove sharp bounds on the enstrophy growth in viscous scalar conservation laws. The upper bound is, up to a prefactor, the enstrophy created by the steepest viscous shock admissible by the $L^\infty$ and total variation bounds and viscosity. This answers a conjecture by D. Ayala and B. Protas (\emph{Physica D}, 2011), based on numerical evidence, for the viscous Burgers equation. 
\end{abstract}

\maketitle

\parskip   2pt plus 0.5pt minus 0.5pt

\section{Introduction}
\label{sec:intro}

We consider the initial-value problem for the one-dimensional viscous Burgers equation
\begin{equation}\label{eq:burgers}
\begin{cases}
\p_t u^\nu+ u^\nu \p_x u^\nu= \nu \p_x^2 u^\nu , & t >0 , \ x \in \mathbb{T} \, ,  \\ 
u^\nu(0,x)=u_0(x) , & x \in \mathbb{T} \, ,
\end{cases}
\end{equation}
where $\nu>0$, $\mathbb{T}=\R / \Z$ is the unit circle equipped with periodic boundary conditions, and $u_0$ has zero average. Solutions to \eqref{eq:burgers} exhibit `steepening of gradients', caused by the advection term, which is subsequently arrested by viscous dissipation when the diffusion begins to dominate. 

We are interested in studying the maximum amplification achieved by the ``enstrophy'' $\mathcal{E}(t) :=\| \p_x u^\nu(t,\cdot) \|_{L^2}^2$. This problem has two variants. Lu and Doering~\cite{MR2482997}, motivated by enstrophy growth in the Navier-Stokes equations and its connection to potential finite-time singularity, proposed to \emph{maximize the instantaneous enstrophy growth rate}, namely,
\begin{equation}
	\label{eq:enstrophyineq}
\frac{1}{2} \frac{\d \mathcal E}{\d t} = - \nu \int (\p_x^2 u^\nu)^2 \, \mathrm d x + \frac{1}{2} \int (\p_x u^\nu)^3 \, \mathrm dx
\end{equation}
%the right-hand side of~\eqref{eq:enstrophyineq},
\emph{subject to the constraint} $\mathcal{E} = \mathcal{E}_0$. The maximizers are explicit and saturate the analytical upper bound
\begin{equation}
\frac{\d \mathcal E}{\d t} \leq \frac{C \mathcal E^{\frac{5}{3}}}{\nu^{1/3}} \, , \quad  \mathcal E_0 \geq 1 \, ,
\end{equation}
thereby resolving the problem.
However, Ayala and Protas~\cite{zbMATH05987827} observed numerically that the growth of Lu and Doering's maximizers is not sustained. They proposed to \emph{maximize $\sup_{t > 0} \mathcal E(t)$ subject to the constraint $\mathcal E(0) = \mathcal  E_0$}. As a proxy, the authors numerically maximize the finite-time enstrophy $\mathcal E(T)$ %subject to $\mathcal E(0) = \mathcal E_0$
and observe the scaling
\begin{equation}
	\label{eq:numerical}
\mathcal E(T) \sim \mathcal E_0^{3/2} \, ,\quad \mathcal{E}_0 \to +\infty \, ,
\end{equation}
for the maximal enstrophy under $O(1)$ viscosity at various times $T$. Subsequently, Pelinovsky exhibited solutions with enstrophy satisfying the scaling~\eqref{eq:numerical} with $T \sim \mathcal{E}_0^{-1/2}$ through the Hopf--Cole transformation (\cite[Theorem 1.1]{MR2988280} and~\cite{1202.2071}).

Notably, the analytical upper bound
\begin{equation}
    \label{eq:badbound}
\sup_{t > 0} \mathcal{E}(t) \leq \frac{C \mathcal{E}_0^3}{\nu^4} \, , \quad \mathcal E_0 \geq 1 \, ,
\end{equation}
deduced in~\cite[Appendix A]{zbMATH05987827}, \emph{does not match the numerical evidence}~\eqref{eq:numerical}. In this paper, we clarify this discrepancy by proving the sharp upper bound.  

It will be convenient to non-dimensionalize time so that either $\nu = 1$ and $\mathcal{E}_0 \to +\infty$, as in~\cite{zbMATH05987827}, or $\mathcal{E}_0 = 1$ and $\nu \to 0^+$, which is the convention we follow. This is accomplished via the rescaling
\begin{equation}
    u_\lambda(t,x) = \lambda u(\lambda t,x) \, , \quad \lambda > 0 \, ,
\end{equation}
with $\lambda = \mathcal{E}_0^{-1/2}$. In our normalization, the analytical upper bound~\eqref{eq:badbound} becomes 
%\begin{equation}
    $\sup_{t>0} \mathcal{E}(t) \leq C \nu^{-4}$,
%\end{equation}
whereas the numerically predicted upper bound $\sup_{t > 0} \mathcal{E}(t) \leq C \mathcal{E}_0^{3/2}$ at unit viscosity becomes
\begin{equation}
    \label{eq:predictedbound}
    \sup_{t>0} \mathcal{E}(t) \leq C \nu^{-1} \, , \quad \nu \in (0,1] \, .
\end{equation}

In \cite[Eq. (2.11)]{MR1840744}, Biryuk proved~\eqref{eq:predictedbound}, but with a constant $C$ depending on higher Sobolev norms of the initial datum. His bound does not apply to the maximization problem in~\cite{zbMATH05987827}.

The upper bound~\eqref{eq:predictedbound} has a natural interpretation. It is the entropy created by the developed viscous shock
\begin{equation}
    u(x) = - U \tanh \frac{x}{\ell}
\end{equation}
with $U = O(1)$ and characteristic length scale $\ell = \nu/U$. In a size $O(\ell)$ neighborhood of the origin, its enstrophy is $O(U^2/\ell) = O(U^3/\nu)$; that is, the viscous shock dissipates energy at the rate $O(U^3)$. Then we anticipate that the most advantageous configuration for enstrophy growth is to develop the steepest viscous shock possible from initial data with $O(1)$ enstrophy. This process is primarily restricted by the total variation $\| \p_x u^\nu(t,\cdot) \|_{L^1}$, which is monotonically decreasing and, on the torus, controlled by the initial enstrophy. We leverage this observation, not exploited in~\cite{zbMATH05987827}, to prove the sharp upper bound.

\begin{theorem}[Upper bound on the enstrophy growth for the viscous Burgers equation]\label{th:main-b}
There exists an absolute constant $C_0 > 0$ such that the following holds. Let $u_0 \in H^1(\T)$ with $\| \p_x u_0 \|_{L^2} \leq 1$ and zero mean. Let $\nu > 0$ and $u^\nu$ be the solution of the viscous Burgers equation \eqref{eq:burgers} on the unit torus $\T$ with initial data $u_0$. Then
\begin{equation}
\| \p_x u^\nu(t,\cdot) \|_{L^2}^2 \leq C_0 (1 + \nu^{-1}) \, , \qquad \text{ for all } t \geq 0 \, .
\end{equation}
\end{theorem}

We further provide an alternative proof (cf. \cite{MR2988280}) of the corresponding lower bound.

\begin{theorem}[Lower bound on the enstrophy growth for the viscous Burgers equation]\label{th:burgers-lowerboundenstrophy}
There exists $u_0 : \T \to \R$ satisfying the conditions of Theorem \ref{th:main-b} and a constant $c_0 > 0$ such that
\begin{equation}
\liminf_{\nu \to 0^+} \nu \| \p_x u^\nu \|_{L^\infty_t L^2_x(\R_+\times \T)}^2 \geq c_0 \, .
\end{equation}
\end{theorem}

Finally, we apply the same strategy as in Theorem~\ref{th:main-b} to prove %the same
sharp upper bounds on the enstrophy growth for solutions to the Cauchy problem for
a multi-dimensional viscous scalar conservation law \begin{equation}\label{eq:scalarconslaw}
\begin{cases}
\p_t u^\nu + \div f(u^\nu) = \nu \Delta u^\nu, & t>0, \ x \in M, \\
u^\nu(0,x) = u_0(x), & x \in M,
\end{cases}
\end{equation}
with locally Lipschitz continuous flux $f : \R \to \R^n$ on the domain $M := \R^n$ or $M := \T_L^n:= (\R/L\Z)^n$ with $L \ge 1$. The key point is to leverage the monotonicity of the $L^\infty$-norm and the total variation. 

\begin{theorem}[Multi-dimensional conservation laws with Lipschitz continuous flux]\label{th:main-multid}
There exists an absolute constant $C_0 > 0$ depending only on the Lipschitz norm of $f|_{[-1,1]}$ and the dimension $n \geq 1$ such that the following holds. Let $M := \R^n$ or $M := \T_L^n$, with $L \geq 1$. Let $u_0 \in L^\infty \cap \dot W^{1,1} \cap \dot H^1(M)$ satisfying
\begin{equation}
\| u_0 \|_{L^\infty} , \  \| \nabla u_0 \|_{L^1} , \ \| \nabla u_0 \|_{L^2} \leq 1 \, .
\end{equation}
Let $\nu > 0$ and $u^\nu$ be the solution of the viscous scalar conservation law~\eqref{eq:scalarconslaw} on $M$ with initial data $u_0$. Then
\begin{equation}
\| \nabla u^\nu(t,\cdot) \|_{L^2}^2 \leq C_0 (1 + \nu^{-1}) \, , \qquad \text{ for all } t \ge 0 \, .
\end{equation}
\end{theorem}
The constant $C_0$ is independent of the domain. We view Theorem~\ref{th:main-b} as an immediate corollary of Theorem~\ref{th:main-multid}.

\section{Proofs}
\label{sec:proofs}

The Cauchy problem~\eqref{eq:scalarconslaw} is globally well-posed in the subcritical space $L^\infty$ in any dimension. The finiteness of the enstrophy, $\dot W^{1,1}$-seminorm, and various other quantities can be proven \emph{a posteriori} by considering the equation~\eqref{eq:equationforgradient} satisfied by $\nabla u^\nu$ with $f'(u)$ viewed as a known bounded function. We refer to~\cite[Theorem 2.9 and Lemma 2.16]{MR1409366} for well-posedness in $L^2 \cap L^\infty$ and $L^\infty \cap H^1$.

Our strategy is based on $L^\infty$ and TV-bounds for the viscous conservation law and two subsequent estimates. First, we propagate the initial enstrophy bound on the time interval $(0,\nu]$; after time $O(\nu)$, the estimate degenerates---see~\eqref{eq:gronwall}. Then we view $u^\nu$ as solving the heat equation with data and right-hand side controlled only by the monotone quantities, and we rely on smoothing, which is effective after time $O(\nu)$, to bound the enstrophy---see~\eqref{eq:firstduhamel} and~\eqref{eq:secondduhamel}. The smoothing estimate depends only on the monotone quantities and controls the solution on $(\nu,+\infty)$. 

We collect the necessary linear estimates in the following lemma. 
\begin{lemma}[Heat estimates]
	\label{lem:heatestimates}
Let $v_0 \in L^\infty \cap \dot W^{1,1}(M)$. Then, for all $t > 0$, we have
\begin{align}
\label{eq:firstestimate}
\| \nabla e^{\nu t \Delta} v_0 \|_{L^2} & \lesssim (\nu t)^{-\frac{1}{4}} \| v_0 \|_{L^\infty}^{\frac{1}{2}} \| \nabla v_0 \|_{L^1}^{\frac{1}{2}} \, , \\ 
\label{eq:secondestimate}
\| \nabla^2 e^{\nu t \Delta} v_0 \|_{L^2} &\lesssim (\nu t)^{-\frac{3}{4}} \| v_0 \|_{L^\infty}^{\frac{1}{2}} \| \nabla v_0 \|_{L^1}^{\frac{1}{2}} \, .
\end{align}
\end{lemma}
All implied constants are allowed to depend on the dimension $n$.

\begin{proof}
By the representation formula via the heat kernel on $\R^n$ (notably, we may consider $L^\infty(M) \subset L^\infty(\R^n)$ in the periodic setting), we have 
\begin{equation}
\| \nabla  e^{\nu t \Delta} v_0 \|_{L^\infty} \lesssim (\nu t)^{-\frac{1}{2}} \| v_0 \|_{L^\infty} \, .
\end{equation}
Since $\nabla$ commutes with $e^{\nu t \Delta}$, we compute 
\begin{equation}
\| \nabla  e^{\nu t \Delta} v_0 \|_{L^1} \lesssim \| \nabla v_0 \|_{L^1} \, .
\end{equation} 
Then the interpolation inequality $\| \cdot \|_{L^2} \leq \| \cdot \|_{L^\infty}^{\sfrac{1}{2}} \| \cdot \|_{L^1}^{\sfrac{1}{2}}$ yields~\eqref{eq:firstestimate}.

To obtain~\eqref{eq:secondestimate}, we combine the $\dot H^1 \to \dot H^2$ smoothing estimate (proven via the Fourier representation and the inequality $|\xi|^s e^{-\nu |\xi|^2 t} \lesssim (\nu t)^{-\sfrac{s}{2}}$, $s \geq 0$)  and~\eqref{eq:firstestimate}:
\begin{equation}
\| \nabla^2 e^{\nu t \Delta} v_0 \|_{L^2} \lesssim (\nu t/2)^{-\frac{1}{2}} \| \nabla e^{\nu t \Delta/2} v_0 \|_{L^2} \overset{\eqref{eq:firstestimate}}{\lesssim} (\nu t/2)^{-\frac{3}{4}} \| v_0 \|_{L^\infty}^{\frac{1}{2}} \| \nabla v_0 \|_{L^1}^{\frac{1}{2}} \, .
\end{equation}
\end{proof}

\begin{proof}[Proof of Theorem \ref{th:main-multid}]
\textbf{Step 1.} \emph{Conserved quantities}. By the maximum principle for the PDE $\p_t u^\nu + f'(u) \cdot \nabla u^\nu = \nu \Delta u^\nu$ and using the viscous continuity equation
\begin{equation}
	\label{eq:equationforgradient}
\p_t \nabla u^\nu + \nabla (f'(u^\nu) \cdot \nabla u^\nu) = \nu \Delta \nabla u^\nu
\end{equation}
for $\nabla u^\nu$, we deduce 
\begin{equation}
\| u^\nu \|_{L^\infty} , \  \| \nabla u^\nu \|_{L^1} \leq 1 \, ,
\end{equation}
respectively; see~\cite[Theorem 2.29]{MR1409366} for further details.

\textbf{Step 2.} \emph{Propagation estimates for the entropy in the time interval $t \in (0,\nu]$}. We multiply~\eqref{eq:equationforgradient} by $\nabla u^\nu$ and integrate by parts to obtain a differential inequality for the enstrophy:
\begin{equation}\label{eq:diffineq2}
\begin{aligned}
\frac{\d}{\d t} \, \frac{1}{2}  \int_{M} |\nabla u^\nu|^2 \, \mathrm d x &= \int_{M} (f'(u^\nu) \cdot \nabla u^\nu) \Delta u^\nu \, \mathrm d x - \nu \int_{M} |\Delta u^\nu|^2 \, \mathrm d x \\
&\leq \frac{\nu^{-1}}{2}  \| f|_{[-1,1]} \|_{{\rm Lip}}^2 \int_{M} |\nabla u^\nu|^2 \, \mathrm d x - \frac{\nu}{2} \int_{M} |\Delta u^\nu|^2 \, \mathrm d x \, ,
\end{aligned}
\end{equation}
where we used Young's inequality and the $L^\infty$-bound %(see \cite[Appendix B.2.d]{MR2597943})
in the second line. By Gronwall's inequality, \eqref{eq:diffineq2} yields
\begin{equation}
\label{eq:gronwall}
\int_{M} |\nabla u^\nu|^2 \, \mathrm d x \leq \| \nabla u_0 \|_{L^2}^2 \exp \left( \nu^{-1} \| f|_{[-1,1]} \|_{{\rm Lip}}^2 \, t \right) \, .
\end{equation}
In particular,
\begin{equation}
	\label{eq:firstineq}
\sup_{t \in (0,\nu]} \| \nabla u^\nu(t,\cdot) \|_{L^2}^2 \lesssim 1 \, .
\end{equation}

\textbf{Step 3.} \emph{Duhamel's formula and smoothing effect: $t \in (\nu, \infty)$}. By Duhamel's formula, we write
\begin{equation}
u^\nu(t_0 + \tau,\cdot) = e^{\nu \Delta \tau} u^\nu(t_0,\cdot) - \int_{t_0}^{t_0+\tau} e^{\nu (t_0+t-s)\Delta} \div f(u^\nu)(s,\cdot) \, \mathrm ds \, , \qquad t_0, \tau \geq 0 \, .
\end{equation}
We bound the right-hand side in terms of the conserved quantities. First, we have
\begin{equation}
	\label{eq:firstduhamel}
\| e^{\nu \Delta \tau} u^\nu( t_0,\cdot) \|_{\dot H^1} \overset{\eqref{eq:firstestimate}}{\lesssim} (\nu \tau)^{-\frac{1}{4}} \| u^\nu( t_0,\cdot) \|_{L^\infty}^{\frac{1}{2}} \| \nabla u^\nu( t_0,\cdot) \|_{L^1}^{\frac{1}{2}} \, .
\end{equation}
Second, we have $\| f(u^\nu) \|_{L^\infty} \leq \| f|_{[-1,1]} \|_{L^\infty}$ and $\| \nabla f(u^\nu) \|_{L^1} \leq \| f|_{[-1,1]} \|_{\rm Lip} \| \nabla u^\nu \|_{L^1}$. Then, moving the divergence onto the heat kernel and applying~\eqref{eq:secondestimate}, we have
\begin{equation}
	\label{eq:secondduhamel}
\begin{aligned}
\left\| \int_{t_0}^{t_0+\tau} \nabla e^{\nu (t_0+\tau-s)\Delta} \cdot f(u^\nu)(s,\cdot) \, \mathrm d s \right\|_{\dot H^1} &\lesssim \int_0^\tau (\nu (\tau-s))^{-\frac{3}{4}} \, \mathrm d s \\ &\lesssim \nu^{-\frac{3}{4}} \tau^{\frac{1}{4}} \, .
\end{aligned}
\end{equation}
We can optimize the above inequalities by choosing $\tau = \nu$. Then
\begin{equation}
	\label{eq:secondineq}
\sup_{t \geq \nu} \| \nabla u^\nu(t,\cdot) \|_{L^2}^2 \lesssim \nu^{-1} \, .
\end{equation}

\textbf{Step 4.} \emph{Conclusion of the proof}.
Combining Steps 2 and 3, we conclude the proof: we use the propagation estimate \eqref{eq:firstineq} to handle times $t  \in (0, \nu]$; and the smoothing estimate \eqref{eq:secondineq} (with $t_0 = t - \nu$) for $t \in (\nu, +\infty)$.
\end{proof}

\begin{remark}[Viscous Burgers equation in 1D]\label{rk:burgers}
In the setting of the Burgers equation \eqref{eq:burgers}, we can avoid using the heat estimates from Lemma~\ref{lem:heatestimates}.  
Indeed, since $\dot W^{1,1}$ and $L^\infty$ have roughly the same ``strength'' in dimension one, we may substitute 
\begin{equation}
\| \p_x e^{\nu t \p_x^2} u(t_0,\cdot) \|_{L^2} \lesssim (\nu t)^{-\frac{1}{4}} \| \p_x u(t_0,\cdot) \|_{L^1} \lesssim (\nu t)^{-\frac{1}{4}}
\end{equation}
for~\eqref{eq:firstduhamel} and
\begin{equation}
\begin{aligned} 
\left\| \int_{t_0}^{t_0+\tau} \p_x  e^{\nu (t-s) \p_x^2} u^\nu \p_x u^\nu \, \mathrm  d s \right\|_{L^2} &\lesssim \int_0^t (\nu (t-s))^{-\frac{3}{4}} \, \mathrm d s \, \| u^\nu \|_{L^\infty_t L^\infty_x} \| \p_x u^\nu \|_{L^\infty_t L^1_x}
\end{aligned}
\end{equation}
for~\eqref{eq:secondduhamel}. On the other hand, the interpolation argument from Lemma~\ref{lem:heatestimates}  is necessary when $n \geq 2$ to maximally utilize the $L^\infty$-bound. %, which is ``stronger'' than $\dot W^{1,1}$ in the sense of scaling. \nncomment{Maybe a hint about the scaling argument?}
\end{remark}

\begin{remark}[Less restrictive assumption]
    From~\eqref{eq:gronwall}, we observe that Theorem~\ref{th:main-multid} remains valid under the assumption $\| u_0 \|_{L^\infty} , \, \| \nabla u_0 \|_{L^1} , \, \nu \| \nabla u_0 \|_{L^2}^2 \leq 1$.
\end{remark}

Finally, we present the proof of Theorem \ref{th:burgers-lowerboundenstrophy}. 

\begin{proof}[Proof of Theorem \ref{th:burgers-lowerboundenstrophy}]
First, we prescribe initial data $v_0$ for the inviscid problem whose entropy solution $v$  shocks only at the origin. Let $v_0 \in C^\infty(\T)$ be odd on the fundamental domain $[-\sfrac{1}{2},\sfrac{1}{2})$, non-negative and concave on $[-\sfrac{1}{2},0)$, increasing on $[-\sfrac{1}{2},-\sfrac{1}{3})$, equal to $+1$ on $[-\sfrac{1}{3},-\sfrac{1}{6}]$, and decreasing on $[-\sfrac{1}{6},0]$. Concavity ensures that the solution is described by the method of characteristics on $[-\sfrac{1}{2},0)$: indeed, given a particle label $\alpha$, the ``local turnover time" $t_*(\alpha) = -1/v_0'(\alpha)$, at which the derivative of the flow map $\eta(\alpha) = \alpha + v_0(\alpha) t$ vanishes, will be at least the time $t_s(\alpha) = v_0(\alpha)/\alpha$ at which the characteristic enters the origin. For $t > t_*(0)$, $v$ has the desired shock at the origin. 

Let $u_0 = U v_0$ where $U = 1/\|\p_x v_0\|_{L^2} > 0$ is a normalizing factor to ensure that the conditions of the above theorem are satisfied. This amounts to a time rescaling: $u(t,x) = U v(Ut,x)$. The whole family of smooth solutions $\{u^\nu\}_{\nu >0}$ to the viscous equation \eqref{eq:burgers} converges to the unique entropy solution $u$ of the inviscid problem
\begin{equation}\label{eq:burgers-inviscid}
\begin{cases}
\p_t u+ u \p_x u= 0  , & t >0  , \ x \in \mathbb{T} \, ,  \\ 
u(0,x)=u_0(x)  , & x \in \mathbb{T} \, .
\end{cases}
\end{equation}
Indeed, compactness in $L^p$, with $p < +\infty$, is ensured by the uniform bounds in $L^\infty$ and $W^{1,1}$ and the convergence along the whole sequence is guaranteed by the uniqueness of the limit entropy solution (owing to Urysohn's subsequence principle); for further details on the compactness argument and on the entropy-admissibility of the limit point, we refer to  \cite[Theorems 4.62, 4.71, and 5.1]{MR1409366}.

Notably, the energy density measure $u^2/2$ satisfies 
\begin{equation}\label{eq:entropydiss}
\p_t \left(\frac{1}{2} u^2 \right) + \p_x \left(\frac{1}{3} u^3\right) = - \lim_{\nu \to 0^+} \nu |\p_x u^\nu|^2 \leq 0 \, .
\end{equation}
The right-hand side of \eqref{eq:entropydiss} is uniformly bounded in $L^1$ and converges in the weak-$\ast$ sense of finite measures. On the time interval $I := (\sfrac{1}{(6U)}+\varepsilon,\sfrac{1}{(3U)}-\varepsilon)$, for $0 < \varepsilon \ll 1$, we have that $u = U(\mathds{1}_{(-U\varepsilon,0)} - \mathds{1}_{(0,U\varepsilon)})$ in the neighborhood $O := (-U\varepsilon,U\varepsilon)$ and, therefore, the left-hand side of \eqref{eq:entropydiss}  is given by $- \sfrac{2}{3}\,  U^3\,  \delta_{\{x=0\}}$ in~$O$.  Hence, we conclude 
\begin{equation}
\liminf_{\nu \to 0^+} \, \nu \| \p_x u^\nu \|_{L^\infty_t L^2_x(I\times O)}^2 \geq \liminf_{\nu \to 0^+} \, \nu \, \frac{1}{|I|} \int_I \int_O |\p_x u^\nu|^2 \, \mathrm dx \, \mathrm dt \geq \frac{2}{3} U^3 \, .
\end{equation}
\end{proof}

\section*{Acknowledgments}

NDN is grateful to G. Fantuzzi for suggesting the problem and engaging in several interesting conversations. He also thanks A. Bressan, G. M. Coclite, E. Marconi, and E. Zuazua for their helpful comments. DA thanks B. Protas for introducing him to this problem.

NDN is a member of the Gruppo Nazionale per l’Analisi Matematica, la Probabilità e le loro Applicazioni (GNAMPA) of the Istituto Nazionale di Alta Matematica (INdAM); he was partially supported by the Alexander von Humboldt-Professorship program and by the Transregio 154 Project ``Mathematical Modelling, Simulation and Optimization Using the Example of Gas Networks'' of the Deutsche Forschungsgemeinschaft. DA was supported by NSF Postdoctoral Fellowship  Grant No.\ 2002023

%\vspace{1cm}
\bibliographystyle{abbrv}
\bibliography{Enstrophy-ref}
\vfill 

\end{document}